# On the History of Number Line

Galina Sinkevich, associate professor, Department of Mathematics,
Saint Petersburg State University of Architecture and Civil Engineering
Vtoraja Krasnoarmejskaja ul. 4, St. Petersburg, 190005, Russia

*Abstract.* The notion of number line was formed in XX c. We consider the generation of this conception in works by M. Stiefel (1544), Galilei (1633), Euler (1748), Lambert (1766), Bolzano (1830-1834), Méray (1869-1872), Cantor (1872), Dedekind (1872), Heine (1872) and Weierstrass (1861-1885).

A number line is an abstract notion which evolved early in 20[th] century. A number line should be distinguished from solid and geometrical line. A solid number line or an interval is an image which came into existence in the ancient world. As a notion, a geometrical line or axis formed in analysis in the period from the 16[th] to 18[th] century. The notion of a right line or a curve as a locus emerged in the 17[th] century in the earliest works on analysis [1].

As a concept, a number line formed in works of Cantor and Dedekind, however, the term itself, first known as a "number scale" and thereafter, as a "number line" has been used since 1912: "Thus, the positive and negative numbers together form a complete scale extending in both directions from zero" [2].

A continuum, the philosophic idea of the continuous or extended, was the prescience of a number curve. It dates back to the ancient world (Zeno, Aristotle), Middle Ages (Boethius), and beginning of the Modern Age (G. Buridan, T. Bradvardin), and thereafter, Leibniz comes.

In the ancient world, numbers were presented as a set of natural numbers. Rational positive proportions of geometric magnitudes were quantities. Irrational quantities $\pi$ and $\sqrt{2}$ were determined through approximants. The reasoning was based on Eudoxus' method of exhaustion. The estimations were assumed as greater and smaller. The approximation techniques reached the peak in works of Archimedes and later, in works of oriental mathematicians. Zero was not regarded as a number for a long time. Although results of certain problems were negative (e.g., Diophantus obtained such negative numbers), such negative numbers were not deemed to be competent; in certain rare cases, they were interpreted as a debt. Before the modern age, only positive routes were sought when solving equations. Irrational numbers conventionally from Euclid were understood as non-extractable radicals.

Irrational number were called (e.g. by Newton) surdi (deaf) or false. Imaginary numbers which first appeared in 1545 in Cardano's works were called sophistic numbers.

Michael Stifel (1487–1567) was the first to define negative numbers as numbers that are less than zero and positive numbers as those that are greater than zero. It was he who described zero, fractional and irrational numbers as numbers. Stifel wrote in what way whole, rational and irrational numbers relate to each other.

Let us address Stifel's book of 1544 "Arithmetica integra" [3]. Stifel recognized that there are infinitely many fractions and irrational numbers between two nearest whole numbers. He considered a unit segment (2, 3) and located infinite sequences therein

$$2\frac{1}{2}, 2\frac{1}{3}, 2\frac{2}{3}, 2\frac{1}{4}, 2\frac{3}{4}, 2\frac{1}{5}, 2\frac{2}{5}, 2\frac{3}{5}, 2\frac{4}{5}, 2\frac{1}{6}, 2\frac{5}{6}, 2\frac{1}{7}, 2\frac{3}{7}, \ldots \text{ and}$$

$$\sqrt{5}, \sqrt{6}, \sqrt{7}, \sqrt{8}, \sqrt[3]{9}, \sqrt[3]{10}, \sqrt[3]{11}, \sqrt[3]{12}, \sqrt[3]{13}, \sqrt[3]{14}, \sqrt[3]{15}, \sqrt[3]{16}, \sqrt[3]{17}, \sqrt[3]{18}, \sqrt[3]{19}, \sqrt[3]{20},$$
$$\sqrt[3]{21}, \sqrt[3]{22}, \sqrt[3]{23}, \sqrt[3]{24}, \sqrt[3]{25}, \sqrt[3]{26}, \sqrt[4]{17}, \sqrt[4]{18}, \sqrt[4]{19}, \sqrt[4]{20}, \sqrt[4]{21}, \sqrt[4]{22}, \sqrt[4]{23}, \sqrt[4]{24}, \sqrt[4]{26}, \ldots$$

[3, p.104]. See [4] for more detail.

In 1596, in keeping with the traditions of German Rechenhaftigkeit (calculability) and following Archimedes whose works won prominence in Europe in the middle of the 16[th] century, Ludolph van Ceulen (1539–1610) calculated a 35-decimal number π. Since then, number π was known as "Ludolphine number" until the end of the 19[th] century.

Unlike Stifel, G. Galilei (1564 – 1642) felt that mathematics and physics are linked with each other; all his reasonings were accompanied by examples from optics, mechanics, etc. A line was understood as a result of movement.

In 1633, in his book entitled "The Discourses and Mathematical Demonstrations Relating to Two New Sciences", Galilei discoursed of the counting distribution: "SIMP. If I should ask further how many squares there are one might reply truly that there are as many as the corresponding number of roots, since every square has its own root and every root its own square, while no square has more than one root and no root more than one square.

But if I inquire how many roots there are, it cannot be denied that there are as many as there are numbers because every number is a root of some square. This being granted we must say that there are as many squares as there are numbers because they are just as numerous as their roots, and all the numbers are roots. Yet at the outset we said there are many more numbers than squares, since the larger portion of them are not squares. Not only so, but the proportionate number of squares diminishes as we pass to larger numbers. Thus up to 100 we have 10 squares, that is, the squares constitute 1/10 part of all the numbers; up to 10000, we find only 1/100 part to be squares; and up to a million only 1/1000 part; on the other hand in an infinite number, if one could conceive of such a thing, he would be forced to admit that there are as many squares as there are numbers all taken together.

So far as I see we can only infer that the totality of all numbers is infinite, that the number of squares is infinite, and that the number of their roots is infinite; neither is the number of squares less than the totality of all numbers, nor the latter greater than the former; and finally the attributes "equal," greater," and "less," are not applicable to infinite, but only to finite, quantities." [5, "First day"].

Numbers considered before the 18th century were natural, rational and irrational numbers (as non- extractable radicals – Kestner [6]), that is to say, algebraic irrational ones. The assumption that number π is irrational was voiced by Arab scientists starting from the 11th century [7, p. 176].

In 1748, speaking of tangents of angles less than $30^0$ in his "Introduction to Analysis of Infinitely Small", Euler used the term "too irrational" (nimis sunt irrationals) [8, p. 120]. In the same work, he made an assumption that in addition to irrational algebraic numbers (i.e. products of algebraic equation with rational coefficients), there are transcendental irrational numbers as well which are obtained as a result of transcendental calculations, e.g. taking logarithms. We would note that symbols π and e established after Euler's "Introduction to Analysis of Infinitely Small" was published, although π could be found in works of various mathematicians starting from the beginning of the 18th century. The symbol π goes back to Greek περι- – around, about; or περι-μετρον – perimeter, circumference; or περι-φέρεια – circumference, arc.

In his work of 1766 which was published in 1770 [9] translated into Russian [10, p. 121–143], Johann Heinrich Lambert (1728–1777) proved irrationality of numbers π and e. His reasoning was in 1794 supplemented by Legendre (Legendre's lemma) [10, p. 145–155, 11]. Fourier studied this issue as well (1815).

In 1767, in its "Sur la resolution des équations numériques" (Berlin, 1769), Lagrange defined irrational numbers as those determined by a nonterminating continued fraction.

In 1821, A. –L. Cauchy defined irrational numbers as limits of converging sequences. However, he did not define procedures or operations therewith [12].

In 1830s, in his manuscript entitled "Theory of Values" (Größenlehre), Bernard Bolzano made an attempt to develop a theory of a real number. Bolzano used the exhaustions method and the notions he stated in 1817 regarding the exact least upper bound, and the sequence convergence criterion [Bolzano, B. Rein analytischer Beweis des Lehrsatzes, daß zwischen zwey Werthen, die ein entgegengesetztes Resultat gewähren, wenigstens eine reelle Wurzel der Gleichung liege. – Prag: Gottlieb Haase. – 1817. – 60 s.][1]. He introduced the notion of a measurable number, relations described as "equal to", "greater than", "less than"; asserted density (pantachisch) of a set of real numbers. Bolzano introduced infinitely great and infinitely small numbers. If his manuscript were published and recognized by contemporaries, we would have probably dealt with another kind of analysis, a nonstandard one. Remember what Putnam said about the emergence of the epsilon-delta language: "If the epsilon-delta methods had not been discovered, then infinitesimals would have been postulated entities (just as 'imaginary' numbers were for a long time). Indeed, this approach to the calculus enlarging the real number system–is just as consistent as the standard approach, as we know today from the work of Abraham Robinson. If the calculus had not been 'justified' Weierstrass style, it would have been 'justified' anyway" [Putnam, 1974]. [13]. The treatment of numbers as variable quantities might have started with Leibnitz' assertion to the effect that $\lim \frac{\Delta y}{\Delta x} = \frac{dy}{dx}$. Along with constant numbers, there are variable numbers in Bolzano's works, both measurable and nonmeasurable. The limit or boundary of such numbers can variable as well. "If a variable, however, measurable number $Y$ constantly remains greater than a variable, however, measurable $X$, and moreover, $Y$ has no least value, and $X$ has no greatest one, then there is at least one measurable number $A$ which constantly lies between the boundaries of $X$ and $Y$.

---

[1] This criterion was called after Cauchy, although in Bolzano's works, it was stated in 1817, while in Cauchy's works it was stated in 1821.

If, further, the difference $Y - X$ cannot infinitely decrease, then there are infinitely many such numbers lying between $X$ and $Y$.

However, if this difference infinitely decreases, then there is one and only one such number. And, finally, if the difference $Y - X$ infinitely decreases and either $X$ has the largest or $Y$ has the least value, then there is no measurable number constantly lying between $X$ and $Y$ [14, p. 525]. There is already a prescience of the notion of a cut here, which is going to appear in 1872 in Dedekind's works. However, in 1830, nobody was aware of the existence of transcendental numbers. Bolzano was more of a philosopher; his ingenious mathematical wisdom was not based on his professional activities, although his philosophical approach to understanding of the continuity formed the line of development of the concept of a number and continuity. The compelled ban from teaching, the lack of academic intercourse and scientific literature prevented his ideas from taking shape of an independent theory. However, his ideas became part of the theory of functions and the theory of sets.

In 1840, G. Liouville started developing the notion and theory of transcendental numbers. In 1844, he published a small article in *Comptes Rendus* where he said that an algebraic number cannot be approximated by a rational fraction [15]. His further research constituted the theory of transcendental numbers. In 1873, Hermite proved the transcendence of number *e* [16]; in 1882, Lindeman proved transcendence of number π [17]; in 1885, his proof was simplified by Weierstrass [18].

However, the theory of a real number had not been created as yet. The terms 'zero', 'greater than', 'less than', or 'equal to zero' could not be rigorously defined. Therefore, in his lectures of 1861 in differential calculus, Weierstrass proved the theorem as follows: "A continuous function where a derivative that is within certain intervals of arguments is always equal to zero amounts to a constant" [19, p. 118].

In 1869, a French mathematician Charles Méray (1835-1911) developed the theory of a real number [20]. Being based on converging sequences and having introduced a relation of equivalence between them, Méray introduced the notion of a nonmeasurable number as a fictive limit: "An invariant converges to a certain fictive nonmeasurable limit if it converges to a point which does not allow for a precise definition. If nonmeasurable limits of two converging variants are equal, these variants will be of equal value" [21, p. 2]. Méray defined the relation of comparison of and operation with nonmeasurable numbers. However, his sophisticated language, clumsy terms; his remoteness from mathematical life of Paris (Méray lived in a village and dealt in winegrowing for many years; thereafter, he gave lectures at Dijon university), his utter antagonism to or lack of knowledge of achievements in mathematics after Lagrange, put restrictions on his theory. He believed that "one will never come across discontinuous functions which have no derivatives and are non-integrable, so no worries about them. There is no point in addressing the Laplace equation, Dirichlet principle, because derivatives are defined, calculated, stirred into differential expressions the way Lagrange wanted them to do, that is, with the help of simple operations" [22]. Therefore his contemporaries did not accept his theory, although now they call the concept of a real number the Méray–Cantor concept.

German mathematicians took the initiative of developing the concept of a real number. In 1872, works of E. Heine "Lectures on the Theory of Functions" [23], G. Cantor "Extension of One Theorem from the Theory of Trigonometric Series" [24, p. 9–17] and R. Dedekind "Continuity and Irrational Numbers" [25] were published.

A professional mathematician and teacher, E. Heine explained the theory of a number in terms of fundamental sequences, having introduced a relation of equivalence and order. His explanation had much in common with Cantor's theory as it was developed in the course of joint discussions with him. He was ahead of his colleagues methodically. Hitherto, the notion of a limit in analysis had often been provided in terms of countable sequences [26, 27].

R. Dedekind approached to the definition of a number as an algebraist, having provided an arithmetic definition of a number. Dedekind considered properties of equality, ordering, density of the set of rational numbers R, (numerical corpus, a term introduced by Dedekind in schedules to Dirichlet's lectures he published). In doing so, he tried to avoid geometric representations. Having defined the relation of 'greater than' and 'less than', Dedekind proved its transitivity, existence of infinitely many other numbers between two different numbers, and that for any number it is possible to make a cut of a set of rational numbers into two classes, so that numbers of either class are less than this number and numbers of another class are greater than this one, in which event the very number which makes such cut may be attributed to as one class or to the other, in which event it would be either the greatest one for the first class or the least for the second one.

Subsequently, Dedekind considered points on a straight line and established the same properties for them as he had just found for rational numbers, thus stating that a point on a straight line corresponds to each rational number. "Each point $p$ of a line separates the line into two parts, so that each point of one part is located to the left of each point of the other part. Now I perceive the essence of continuity in the opposite principle, i.e. it is as follows: "If all points of a line fall into such two classes that each point of the first class lies left of each point of the second class, then there is one and only one point which separates the line into two classes, and this is a cut of the line into two fragments. If the system of all real numbers is separated into two such classes that each number of the first class is less than each number of the second class, then there is one and only one number making such cut" [25, p. 17–18].

Dedekind called this property of a line an axiom accepting which we make the line continuous. In this case, Dedekind asserted that it was our mental act which occurred whether the real space is continuous or discontinuous, that such mental completion by new points did not affect the real existence of the space.

The Dedekind's definition of a number is still used in courses of analysis as logically and categorically impeccable. However, as Cantor noted, one cannot use the notion of a number as a cut in analysis. As soon as it concerns an irregular set, this definition is useless.

Having introduced the notion of a number based on fundamental countable sequences just as Heine, G. Cantor passed on Heine: he defined the notion of a limiting point and introduced the hierarchy of limiting sets. In his work of 1874, in his article entitled "Über eine Eigenschaft des Inbegriffes aller reellen algebraiscen Zahlen" [24, p. 18–21], he proved the countability of a set of algebraic irrational numbers and uncountability of a set of real numbers and, therefore, transcendental numbers. Set-theoretic terminology had not formed by that time, the notion of countability would appear in his works later, so he spoke of one-to-one correspondence and used the term 'Inbegriff' (totality) instead of a set. Cantor postulated the one-to-one correspondence between numbers and points on a line and asserted that this could not be proved.

By 1878, Cantor switched from analysis of point spaces to the notion of potency (power of set), formed a continuum hypothesis, considered continuous mapping between sets of various dimensionalities. The more he felt the insufficiency of the definition of continuity through the cut. In his third article of 1878 entitled "Ein Beitrag zur Mannigfaltigkeitslehre" (To the Theory of Manifolds) [24, p. 22–35], he already provided the notions of potency and one-to-one correspondence between manifolds of various dimension. It was in the same article that he introduced the notion of a 'second potency', which means that the continuum hypothesis started to form. See [28] for more detail. According to Cantor, a continuum of numbers is a perfect connected well-ordered set.